\documentclass[a4paper,10pt]{amsart}
\usepackage{amsmath,amssymb,amsthm,latexsym}
\usepackage{verbatim,enumerate}
\usepackage{mdwlist}
\usepackage[dvips]{graphicx}
\usepackage{psfrag}
\usepackage{array}

\usepackage{changebar}

\newtheorem{thm}{Theorem}[section]
\newtheorem{lem}[thm]{Lemma}
\newtheorem{cor}[thm]{Corollary}
\newtheorem{prop}[thm]{Proposition}
\newtheorem{Rem}[thm]{Remark}

\newtheorem{DEF}[thm]{Definition}
\theoremstyle{remark}
\newtheorem*{rem}{Remark}
\newtheorem*{rems}{Remarks}

\newcommand{\rco}{{recursive coatom ordering }}
\newcommand{\rcok}{{recursive coatom ordering}}

\newcommand{\pet}{\vspace{5pt}}

\newcommand{\tn}{\textnormal}
\newcommand{\tp}{\hat{1}}
\newcommand{\bm}{\hat{0}}

\newcommand{\Bpi}{\textnormal{Bier}(P,I)}
\newcommand{\Bli}{\textnormal{Bier}(L,I)}
\newcommand{\BBli}{{\overline{\textnormal{Bier}(L,I)}}}
\newcommand{\BBpi}{{\overline{\textnormal{Bier}(P,I)}}}
\newcommand{\Blil}{{\textnormal{Bier}(L,I)_{<\hat{1}}}}
\newcommand{\nin}{\noindent}
\newcommand{\cf}{{\mathcal F}}
\newcommand{\cs}{{\mathcal S}}
\newcommand{\cfh}{{\hat{\mathcal F}}}
\newcommand{\cfa}{\hat{{\mathcal F}_\alpha}}
\newcommand{\co}{{\mathcal O}}
\newcommand{\coa}{{\mathcal O}_\alpha}
\newcommand{\cl}{{\mathcal L}}
\newcommand{\cg}{{\mathcal G}}
\newcommand{\cn}{{\mathcal N}}
\newcommand{\cc}{{\mathcal C}}
\newcommand{\ca}{{\mathcal A}}
\newcommand{\bl}{{\rm{Bl}\,}}
\newcommand{\bla}{{\rm{Bl}\,}_\alpha}
\newcommand{\blb}{{\rm{Bl}\,}_B}
\newcommand{\Iff}{ if and only if }
\newcommand{\coat}{\textrm{coat}(\cfa)}
\newcommand{\bdm}{\begin{displaymath}}
\newcommand{\edm}{\end{displaymath}}
\newcommand{\nsub}{\nsubseteq}
\newcommand{\nsup}{\nsupseteq}
\newcommand{\sd}{{\rm sd}}
\newcommand{\mf}{{\rm mf}}

\title{Simplicial shellable spheres via combinatorial blowups}

\author{Sonja Lj. \v{C}uki\'{c}}
\address{Department of Computer Science, Eidgen\"ossische Technische
Hochschule, Z\"urich, Switzerland}
\email{sonja.cukic@inf.ethz.ch}
\author{Emanuele Delucchi}
\address{Department of Mathematics, Eidgen\"ossische Technische
Hochschule, Z\"urich, Switzerland}
\email{delucchi@math.ethz.ch}
\thanks {Research partially supported by TH-Projekt 0-20268-05, and by the Swiss National Science
  Foundation, project PP002--106403/1}

\date{\today}

\begin{document}
\begin{abstract}
The construction of the Bier sphere $\textrm{Bier}(K)$ for a
simplicial complex $K$ is due to Bier (\cite{B}, \cite{M}). Bj\"orner,
Paffenholz, Sj\"ostrand and Ziegler \cite{BPSZ} generalize this
construction to obtain a Bier poset $\textrm{Bier}(P,I)$ from any
bounded poset $P$ and any proper ideal $I\subseteq P$. They show shellability
of $\textrm{Bier}(P,I)$ for the case $P=B_n$, the boolean lattice, and obtain thereby 'many shellable spheres' in the sense of Kalai \cite{Ka}. 

We put the Bier construction into the general framework of the theory
of nested set complexes of Feichtner and Kozlov \cite{FK}. We obtain
'more shellable spheres' by proving the general statement that
combinatorial blowups, hence stellar subdivisions, preserve shellability.
\end{abstract}

\maketitle

\section*{Introduction}

Let $K$ denote an abstract simplicial complex on the vertex set $[n]:=\{1,\dots ,n\}$, and write $\mathcal{F}(K)$ for the poset of its faces. The (combinatorial) {\em Alexander dual} for $K$ is the simplicial complex $A(K)$ whose faces are the complements (in $[n]$) of the {\em non-}faces of $K$. Thus, $A(K)=\{[n]\setminus \sigma \mid \sigma \in 2^{[n]} \setminus\, K\}$. The topological motivation for this suggestive name is that, in fact, $K$ and $A(K)$ can be 'put together' to a sphere. A very nice construction of this sphere is due to Thomas Bier \cite{B} and can be found in \cite[p. 111-116]{M}. Bier showed that the deleted join $\textrm{Bier}_n(K):=(K\ast A(K))_\Delta$ is an $(n-2)$-sphere with at most $2n$ vertices \cite[Theorem 5.6.2]{M}. The idea behind this proof is to embed $\mathcal{F}(K)$ in the boolean lattice $B_n$ and see that $\textrm{Bier}_n(K)$ is in fact a subdivision of the boundary of the $(n-1)$-simplex $\Delta^{(n-1)}=\Delta(B_n)$.

Bj\"orner et al.\ generalize this construction in \cite{BPSZ} to obtain a {\em Bier
  poset} $\textrm{Bier}(I,P)$ associated to any proper lower ideal $I$ in any
bounded poset $P$. They show that, for any such $P$, the order complex
  of $\overline{\textrm{Bier}(P,I)}$ is PL-homeomorphic to that of
  $\overline{P}$. In the boolean case $P=B_n$ we have
  $\textrm{Bier}(B_n, I)=\textrm{Bier}_n(I)$ (at the right hand side of the equality $I$ is seen as an abstract simplicial complex). In the same paper
  shellability of $\textrm{Bier}(B_n,I)$ is proven, together with a
  characterization of its $g$-vector. In particular, it is pointed out
  that, for large $n$, this construction leads to 'many simplicial shellable $(n-2)$-spheres', most of them lacking convex realization (see \cite{Ka}).

We put this construction in the context of the theory of nested set
complexes (developed by Feichtner and Kozlov in \cite{FK}). Specifically, we find a conceptual way of proving that $P$ and
$\textrm{Bier}(P,I)$ are PL-homeomorphic if $P$ is a lattice. Moreover, we obtain shellability of $\textrm{Bier}(P,I)$ for {any} shellable lattice $P$ and 'even more' simplicial shellable spheres.

In section \ref{defs} we introduce notations and define the basic
characters of this paper. Section \ref{building} relates the poset $P$
to a building set in the poset $\textrm{Bier}(P,I)$ and then,
restricting to the case when $\textrm{Bier}(P,I)$ is a semilattice,
shows how combinatorial blowups relate order complexes of $\overline{P}$ and $\overline{\textrm{Bier}(P,I)}$. 
 The core of Section
\ref{shell} is Theorem \ref{thm_rco}, where we show that combinatorial
blowups preserve shellability. This applies in particular to the Bier
construction and implies that $\textrm{Bier}(P,I)$ is shellable and 
homotopy equivalent to a wedge of spheres whenever $P$ is a shellable poset. Note that this allows even to iterate the Bier construction to get shellable triangulations of spheres with any number of vertices, whose isomorphism type depends on $P$, $I$, and the number of iterations.

{\bf Acknowledgments.} We would like to thank Eva Maria Feichtner and Dmitry Kozlov for pointing out this problem to us, and for their useful suggestions and comments.

\section{Notations, definitions and basic properties.}\label{defs}

\subsection{Posets.}$\,$

\pet

\nin In this section we will give a summary of the standard definitions and notations which will be used further in the paper. For a general reference to the theory of posets and lattices, we refer the reader to \cite[Chapter 3]{Sta}.

Let $(P,\leq)$ be a poset. All posets considered in this paper will be
 finite. We call $P$ \emph{bounded} if there exist elements $\bm, \tp
 \in P$ so that $\bm\leq x \leq \tp$ for all $x\in P$. We will write
 $\overline{P}$ for the proper part of $P$, that is
 $\overline{P}=P\setminus \{\bm,\tp\}$. Also, let $P_{\leq x}=\{y\in
 P\mid y\leq x\}$. Similarly, for any $\cg\subseteq P$ with order induced by $P$, let us write $\cg_{\leq x}=\{y\in
 P\mid y\in \cg,\ y\leq x\}$.  We say that $y$ \emph{covers} $x$ if
 $y>x$ and there is no $z\in P$ so that $y>z>x$; in this case we will
 write $x\lessdot y$. For $x\leq y$, the \emph{interval} $\{z\in P\mid
 x\leq z\leq y\}$ is denoted by $[x,y]$. A subset $I\subseteq P$ is an \emph{ideal} in $P$ if, for all $y\in I$ and $x\in P$, if $x\leq y$, then $x\in I$. An ideal $I$ is called \emph{proper} if $I\neq P$ and $I\neq \emptyset$. With $\Delta(P)$ we will denote the \emph{order complex} of the given poset $P$: the abstract simplicial complex whose vertices are the elements of $P$, and faces are all chains in $P$ (including the empty chain). In this paper we will assume that the empty face is an element of every non-empty abstract simplicial complex.
  
A poset $\cl$ is called a \emph{meet-semilattice}, or simply
\emph{semilattice}, if every pair of elements $x,y\in \cl$ has a
unique maximal lower bound, which is called \emph{meet} of these two
elements, and is denoted by $x\wedge y$. All semilattices have a unique minimal element called $\bm$, and for any $A=\{a_1,\cdots,a_t\}\subseteq \cl$, the set $\{x\in \cl\mid x\geq a_i,\tn{ for all }i\in [t]\}$ is either empty, or it has a unique minimal element, its \emph{join}, $\bigvee A=a_1\vee\dots\vee a_t$.  Finally, 
a semilattice $L$ is a {\em lattice} if meet and join are defined for any pair of elements of $L$.

\begin{DEF}\tn{(\cite[Definition 1.1]{BPSZ})} Let $P$ be a bounded poset of finite length and $I\subset P$ a proper ideal. Then the poset $\Bpi$ is defined as follows:
\begin{itemize}
\item elements are all intervals $[x,y]\subseteq P$ such that $x\in I$ and $y\notin I$, together with an additional top element $\tp$,
\item intervals are ordered by reverse inclusion, i.e. $[x,y]\leq [v,w]$ if and only if $x\leq v<w\leq y$.
\end{itemize}
\end{DEF} 

\begin{figure}[ht]
\begin{center}
\psfrag{P}{${P}$}
\psfrag{I}{${I}$}
\psfrag{y}{$\scriptstyle{y}$}
\psfrag{x}{$\scriptstyle{x}$}
\psfrag{v}{$\scriptstyle{v}$}
\psfrag{w}{$\scriptstyle{w}$}
\psfrag{0}{$\scriptstyle{\bm}$}
\psfrag{1}{$\scriptstyle{\tp}$}

\includegraphics[scale=0.9]{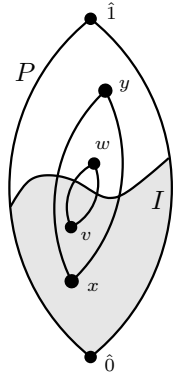}
\caption{Two elements $[x,y]$ and $[v,w]$ of $\Bpi$ with $[x,y]\leq [v,w]$.}
\end{center}
\end{figure}

Clearly, $\Bpi$ is a bounded poset. Furthermore, by \cite[Lemma 1.2]{BPSZ}, if $P$ is a lattice, then $\Bpi$ is also a lattice.
\newpage
\subsection{Building sets, nested sets, combinatorial blowups and stellar subdivisions.} $\,$
\pet

The theory of building sets, nested set complexes and combinatorial
blowups for general semilattices was initiated and developed by Feichtner and Kozlov, \cite{FK}, as the combinatorial framework of the {\em wonderful models for subspace arrangements} by de Concini and Procesi \cite{dCP, F1}. However, this theory has found application in many different contexts, as can be seen in \cite{F2,FK2,FM,FS}. We recall here the basic definitions and refer to \cite{FK} for a comprehensive introduction.
\begin{DEF}\label{def_building} \tn{(\cite[Definition 2.2]{FK})} Let $\cl$ be a semilattice. A subset $\cg$ of~$\cl{\setminus}\{\bm\}$  is called a {\bf building set} of $\cl$ if for any $x\in \cl\setminus \{\bm\}$ and $\max \cg_{\leq x}=\{x_1,\dots,x_t\}$, there is an isomorphism of posets
\begin{equation*}\psi_x:\prod_{i=1}^t[\bm,x_i] \longrightarrow [\bm,x],
\end{equation*}
so that $\psi_x(\bm,\dots,\bm,x_i,\bm,\dots,\bm)=x_i$, for all $i\in [t]$.
\end{DEF}

\begin{rem} Note that this definition does not really require the semilattice structure of $\cl$. Therefore such an object can be defined in any bounded poset with a unique minimal element. However, if $\cl$ is a semilattice, then $\psi_x$ can always be chosen to be the canonical map $(x_1,\dots,x_t) \mapsto x_1 \vee \dots \vee x_t $.
\end{rem}

\begin{DEF} \label{def_nested} \tn{(\cite[Definition 2.7]{FK})} Let $\cl$ be a semilattice and $\cg$ a building set of $\cl$. A (possibly empty) subset $N$ of $\cg$ is called {\bf nested} if for any $\{x_1,\dots,x_t\}\subseteq N$, where $t\geq 2$ and any two distinct elements $x_i$ and $x_j$ are incomparable, the join $x_1\vee\dots\vee x_t$ exists and does not belong to $\cg$. 

The  nested sets in $\cg$ form an abstract simplicial complex, called the
{\bf nested set complex} of $\cg$ in $\cl$, and which will be denoted by  $\cn(\cl,\cg)$. 
\end{DEF}
If it is clear which semilattice $\cl$ is meant, we will write $\cn(\cg)$ instead of $\cn(\cl,\cg)$.
\begin{rem}
It is not hard to see that if $\cg$ is the maximal building set in the given semilattice $\cl$, then $\cn(\cl,\cg)=\Delta({\cl}\setminus \{\bm\})$.
\end{rem}
\begin{DEF}\label{def_Bl} \tn{(\cite[Definition 3.1]{FK})} For a semilattice $\cl$ and an element $\alpha\in \cl$ we define a new poset $\bl_\alpha \cl$ on the set of elements $$ \{x\in \cl\mid x\ngeq \alpha\}\cup\{[\alpha,x]\mid x\in \cl,\ x\ngeq \alpha\tn{ and }(x\vee \alpha)_{\cl} \tn{ exists}\},$$ with order relation defined as follows:
\begin{enumerate}[(i)]
  \item $y>z$ in $\bl_\alpha \cl$ if $y>z$ in $\cl$;
  \item $[\alpha,y]>[\alpha,z]$ in $\bl_\alpha \cl$ if $y>z$ in $\cl$;
  \item $[\alpha,y]>z$ in $\bl_\alpha \cl$ if $y\geq z$ in $\cl$;
  \end{enumerate}
where in all three cases $y,z\not\geq\alpha$. 

The poset $\bl_\alpha \cl$ is called the {\bf combinatorial blowup of $\cl$ at $\alpha$}.
\end{DEF}
\begin{figure}[ht]
\begin{center}
\psfrag{x}{$\scriptstyle{x}$}
\psfrag{y}{$\scriptstyle{y}$}
\psfrag{z}{$\scriptstyle{z}$}
\psfrag{xy}{$\scriptstyle{xy}$}
\psfrag{xz}{$ \scriptstyle{xz}$}
\psfrag{yz}{$\scriptstyle{yz=\alpha}$}
\psfrag{xyz}{$ \scriptstyle{xyz}$}
\psfrag{0}{$ \scriptstyle{\bm}$}
\psfrag{a}{$ \scriptstyle{[\alpha,\bm]}$}
\psfrag{ax}{$ \scriptstyle{[\alpha,x]}$}
\psfrag{ay}{$ \scriptstyle{[\alpha,y]}$}
\psfrag{az}{$ \scriptstyle{[\alpha,z]}$}
\psfrag{axy}{$ \scriptstyle{[\alpha,{xy}]}$}
\psfrag{axz}{$ \scriptstyle{[\alpha,{xz}]}$}
\psfrag{Bl}{$\bl_\alpha$}
\includegraphics[scale=1]{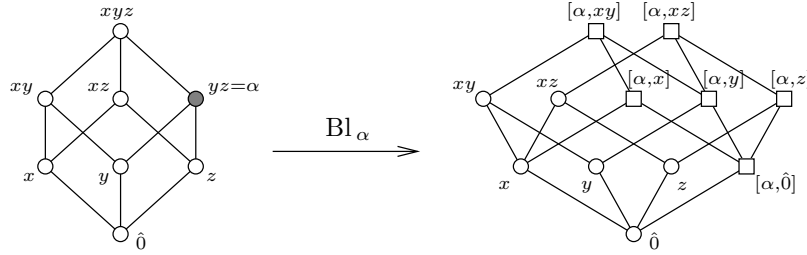}
\caption{Example of the combinatorial blowup of boolean lattice $B_3$ at the element $\{y,z\}$.}
\end{center}
\end{figure}
\begin{rem}
Given a semilattice $\cl$ and an element $\alpha\in \cl$, the poset $\bla\cl$ is again a semilattice, see \cite[Lemma 3.2]{FK}.
\end{rem}
{\bf Example.} Let $\cf=\cf(K)$ be the face semilattice of some simplicial complex $K$, that is, elements of $\cf(K)$ are faces of $K$, and they are ordered by inclusion (the minimal element of $\cf(K)$ is the empty face of $K$). Let $\alpha\in \cf$ be a face of $K$. By the previous remark, $\bla\cf$ is again a  semilattice. If we have two elements from  $\bla\cf$, let us see what their meet is equal to:

$\bullet$ If $F,G\in \cf\cap\bla\cf$, then  $F\wedge G=F\cap G$, seen as an element of $\bla\cf$.

$\bullet$ If $F,[\alpha,G]\in \bla\cf$, where $F,G\in \cf$, then $F\wedge [\alpha,G]=F\cap G\in \bla\cf$.

$\bullet$ Finally, if $[\alpha,F],[\alpha,G]\in \bla\cf$, $F,G\in \cf$, then  $[\alpha,F]\wedge[\alpha,G]=[\alpha,F\cap G]$.\\

We proceed with the definition of stellar subdivision for abstract simplicial complexes. Note that, passing to the geometric realization, this translates to the well-known corresponding geometrical notion.

\begin{DEF}\label{def_stellar}
The {\bf stellar subdivision} of  a simplicial complex $K$ with respect to a non-empty face $F$ is the simplicial complex $\sd_F(K)$ whose faces are $$\{G\in K\mid G\nsup F\}\cup \{G\cup\{v_F\}\mid G\in K,\ G\nsup F,\tn{ and }G\cup F\in K\}.$$
\end{DEF}

\begin{rems}$ $

 $\bullet$ It was noticed in \cite[Section 3]{Ko} that $\cf(\sd_F(K))=\bl_F(\cf(K))$, that is, stellar subdivisions are instances of combinatorial blowups.
 
 $\bullet$ It is known that there exists a sequence of elementary collapses and elementary expansions leading from a simplicial complex $K$ to the complex $\sd_F(K)$. In other words, $K$ and $\sd_F(K)$ have the same simple homotopy type, see for example  \cite[Section 3]{Ko} for description of formal deformation from $K$ to $\sd_F(K)$.
 \end{rems}

\section{Building sets in Bier lattices}\label{building}

From now on, unless stated otherwise, we will assume that $L$ is a lattice.
Then $\Bli$ is also a lattice, and we can therefore apply the theory of nested set complexes. We begin by describing a building set in $\Blil$ that is naturally associated to $L$.

\begin{prop}\label{prop_buildingBier}
For any Bier lattice $\Bli$, where $I\subseteq L$ is a proper ideal, 
\bdm
\cg:=\{[\bm,y]\mid y\in \overline{L}\setminus I\}\cup\{[x,\tp]\mid x\in I\setminus \{\bm\}\}
\edm
is a building set in $\Blil$.

\begin{figure}[ht]
\begin{center}
\psfrag{P}{${P}$}
\psfrag{I}{${I}$}
\psfrag{y}{$\scriptstyle{y}$}
\psfrag{x}{$\scriptstyle{x}$}
\psfrag{0}{$\scriptstyle{\bm}$}
\psfrag{1}{$\scriptstyle{\tp}$}
\psfrag{a}{\textnormal{(a)}}
\psfrag{b}{\textnormal{(b)}}
\psfrag{i}{\it (i)}
\psfrag{ii}{\it (ii)}
\psfrag{iii}{\it (iii)}
\includegraphics[scale=0.9]{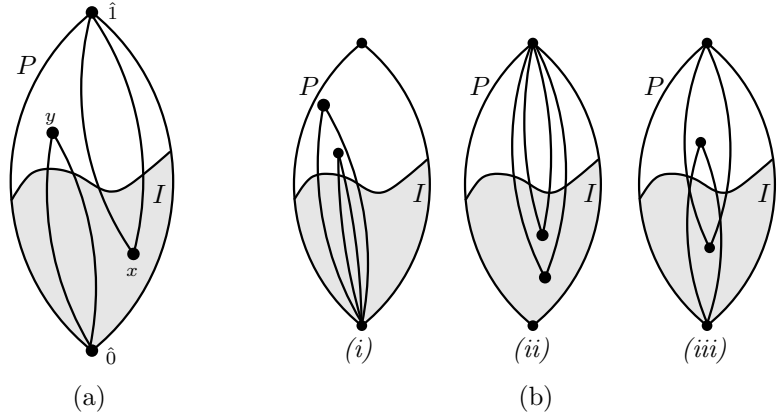}
\caption{(a) Elements of the building set $\cg$. (b) The three cases considered in Lemma \ref{lem_nested}.}
\end{center}
\end{figure}

\begin{proof} Let $[x,y]\in \BBli$. If $[x,y]\in \cg$, then Definition \ref{def_building} is clearly satisfied. Suppose then that $[x,y]\notin \cg$, i.e. $x\neq \bm$ and $y\neq \tp$. It is clear that \
\bdm \max \cg_{\leq [x,y]}=\{[x,\tp],[\bm,y]\}.\edm

Since in $\Bli$ we have that $\left[[\bm,\tp],[x,\tp]\right]=\{[x',\tp]\mid x'\leq x\}$, $\left[ [\bm,\tp],[\bm,y]\right] =\{[\bm,y']\mid y'\geq y\}$, and $\left[ [\bm,\tp],[x,y]\right] =\{[x',y']\mid x'\leq x,\ y'\geq y\}$, it is easy to see that the function 
\begin{eqnarray*}
 \psi_{[x,y]}:\left[[\bm,\tp],[x,\tp]\right]\times\left[[\bm,\tp],[\bm,y]\right]&\longrightarrow& \left[[\bm,\tp],[x,y]\right], \tn{ where}\\
 ([x',\tp],[\bm,y'])&\longmapsto&[x',y'],
\end{eqnarray*}
is an order-preserving bijection of these posets. 

It is also easy to see that
$\psi_{[x,y]}([\bm,\tp],[\bm,y'])=[\bm,y']$ and
$\psi_{[x,y]}([x',\tp],[\bm,\tp])=[x',\tp]$. Therefore the function
$\psi_{[x,y]}$ satisfies the conditions of Definition \ref{def_building}, and $\cg$ is a building set in $\Blil$.
\end{proof}
\end{prop}
It is also true that $\cg$ is a building set in $\Blil$ if $L$ is any
bounded poset, since the definition we gave is independent of the existence of meet and join operations.
However, in that case it is not clear how to characterize the concept of nested sets, even for the special case of Bier posets.

In the lattice case the standard definition works, and therefore we move towards our next goal, the characterization of the nested set complex of $\cg$, which will be reached in Proposition \ref{nestedsetcomplex}. We need a preparatory lemma.

\begin{lem} \label{lem_nested} Let $\cg$ be the building set of $\Blil$ defined in Proposition \ref{prop_buildingBier}. Then a (possibly empty) subset $A$ of $\cg$ is nested if and only if the following three conditions are satisfied:
\begin{enumerate}[(i)]
\item If $[\bm,y_1],[\bm,y_2]\in A$, then $y_1$ and $y_2$ are comparable in $L$.
\item If $[x_1,\tp],[x_2,\tp]\in A$, then $x_1$ and $x_2$ are comparable in $L$.
\item If $[x,\tp],[\bm,y]\in A$, then $x<y$ in $L$.
\end{enumerate}
\begin{proof}
$(\Longrightarrow)$ Suppose that $A$ is nested, and suppose that (i) is not true, i.e. that there exists a pair of elements $[\bm,y_1],[\bm,y_2]\in A$ so that $y_1$ and $y_2$ are not comparable. Then $\left\{[\bm,y_1],[\bm,y_2]\right\}$ is a set of incomparable elements in $A$. But
\bdm
[\bm,y_1]\vee [\bm,y_2]=\left\{\begin{array}{ll} [\bm,y_1\wedge y_2]\in \cg, & \tn{if }y_1\wedge y_2\notin I;\\
\tn{doesn't exist},& \tn{otherwise.}\end{array} \right.
\edm
This is a contradiction with $A$ being nested. Condition (ii) is proved analogously. 

It is left to prove that (iii) is true. Since $[x,\tp]$ and $[\bm,y]$ are incomparable, $[x,\tp]\vee [\bm,y]=[p,q]\in \Blil\setminus \cg$, where $x\leq p<q<\tp$ and $\bm< p<q\leq y$. Therefore $x<y$.

$(\Longleftarrow)$ Suppose that
$A=\{[x_1,\tp],\dots,[x_l,\tp],[\bm,y_1],\dots,[\bm,y_k]\}$, and
$x_1<x_2<\cdots<x_l<y_1<\cdots<y_k$. If $l=0$ or $k=0$, then $A$ is a
chain, and hence nested. Otherwise, it is clear that the cardinality of any set of incomparable elements is at most two. Take any two incomparable elements $[x_p,\tp]$ and $[\bm,y_q]$ from $A$. Then $[x_p,\tp]\vee[\bm,y_q]=[x_p,y_q]\notin \cg$. Therefore $A$ is a nested set.
\end{proof}
\end{lem}

We are now ready to prove the following:

\begin{prop}\label{nestedsetcomplex}
Let $\cg$ be the building set of $\Blil$ defined in Proposition \ref{prop_buildingBier}. Then $\cn(\cg)=\Delta(\overline{L})$.
\begin{proof}
From the proof of Lemma \ref{lem_nested}, we know that if $A\in \cn(\cg)$, then $A=\{[x_1,\tp],\dots,[x_l,\tp],[\bm,y_1],\dots,[\bm,y_k]\}$, for some $x_1<x_2<\cdots<x_l<y_1<\cdots<y_k$ with $x_i\in I\setminus\{\bm\}$, for $i \in [l]$, and $y_j\in \overline{L}\setminus I$, for $j\in [k]$ (where one of $k$ and $l$ can be zero). Define a function $f:\cn(\cg)\to\Delta(\overline{L})$ in the following way: $$f(A):=\{x_1<x_2<\cdots<x_l<y_1<\cdots<y_k\}.$$
The function $f$ is clearly well-defined, injective, and order-preserving. To prove surjectivity, let $S=\{z_1<\cdots<z_s\}\in \Delta(\overline{L})$ and $i:=\max (\{j\mid z_j\in I \}\cup \{0\})$. Then it is easy to see that, since $I$ is an ideal, $S_f=\{[z_1,\tp],\dots,[z_i,\tp],[\bm,z_{i+1}],\dots,[\bm,z_s]\}\subseteq \cg$,  $S_f\in \cn(\cg)$ by Lemma \ref{lem_nested}, and $f(S_f)=S$. 
 
We conclude that $f$ is an order-preserving bijection.
\end{proof}
\end{prop}

The next proposition is similar in spirit to Proposition 4.2 of \cite{FM}, but works in the abstract case as well and does not assume atomicity of the lattice. It describes the behavior of nested set complexes under extension of the building set.

\begin{prop}\label{prop_building_blowup}
Let $\cl$ be a semilattice, and let $\cg$ be a building set in $\cl$. If $\alpha\in \max (\cl\setminus \cg)$, and $\cg'=\cg\cup\{\alpha\}$, then $\cf(\cn(\cg'))=\bl_B(\cf(\cn(\cg)))$, where $B=\max \cg_{\leq \alpha}$.
\begin{proof}
Note first that the number of elements of $B$ is at least $2$, since
otherwise Definition \ref{def_building} would not be satisfied for the building set
$\cg$ and the element $\alpha\in \cl$. If $\{\beta_1,\dots,\beta_t\}\subseteq B$ is an antichain, with $t\geq 2$, then $\beta_1<\vee_{i=1}^t \beta_i\leq \alpha$ and hence $\vee_{i=1}^t \beta_i \in \cl\setminus \cg$. Therefore $B\in \cn(\cg)$. Let us now prove that $\cg'$ is a building set. Since $\alpha$ is a maximal element of $\cl\setminus \cg$, it is easy to see that, for all $x\in \cl\setminus \cg'$, $\max \cg'_{\leq x}=\max \cg_{\leq x}$, and for $x\in \cg'$, $\max \cg'_{\leq x}=x$. Since $\cg$ is  a building set, by Definition \ref{def_building}, $\cg'$ is also a building set.

Define now a map $f:\cf(\cn(\cg'))\to \bl_B(\cf(\cn(\cg)))$ in the following way: 
\bdm
f(A)=\left\{\begin{array}{ll} A, & \tn{if }\alpha \notin A; \\ \left[ B,A\setminus\{\alpha\}\right], &\tn{otherwise,} \end{array}\right.
\edm
where $A\in \cf(\cn(\cg'))$. Let us prove that $f$ is an order-preserving bijection.

$\bullet$ The map $f$ is well-defined:

(1) If $\alpha\notin A$, then let $\{a_1,\dots,a_t\}$, $t\geq 2$, be a set of incomparable elements in $A$, if such exists. Then $\vee_{i=1}^t a_i\in \cl\setminus \cg'\subset \cl\setminus \cg$, and hence $A\in \cn(\cg)$. Since $B$ is an antichain with at least two elements, $\bigvee B=\alpha$, and by assumption $A\in\cn(\cg')$, it follows that $A\nsup B$. Therefore $f(A)=A\in \blb (\cf(\cn(\cg)))$.

(2) If $\alpha \in A$, then $A\setminus\{\alpha\}\in \cn(\cg')$, and therefore $A\setminus\{\alpha\}\in \cn(\cg)$ and $B\nsub A\setminus \{\alpha\}$. In order to prove that $[B,A\setminus\{\alpha\}]$ is an element of $\blb(\cf(\cn(\cg)))$, we need to check that $B\vee (A\setminus \{\alpha\})\in \cf(\cn(\cg))$, i.e. that $B\cup  (A\setminus \{\alpha\})$ is a nested set in $\cg$. 

Note that, for all $x\in A\setminus \{\alpha\}$, $x$ has to be comparable with $\alpha$, since otherwise $\{x,\alpha\}$ would be an antichain in $A$, and $x\vee \alpha$ either does not exist in $\cg'$, or $x\vee \alpha>\alpha$, and hence $x\vee \alpha\in \cg'$. This would contradict the fact that $A$ is a nested set in $\cg'$.

Let $\{x_1,\cdots,x_t\}$, where $t\geq 2$, be a set of incomparable elements in $B\cup  (A\setminus \{\alpha\})$. If there exists some $i\in [t]$ so that $x_i\geq \alpha$, then it is easy to see that $\{x_1,\cdots,x_t\}\cap B=\emptyset$,  and hence $\{x_1,\cdots,x_t\}\subseteq A\setminus \{\alpha\}$. Then $\vee_{i=1}^t x_i\in \cl\setminus \cg$ follows from the fact that $A\setminus\{\alpha\}$ is nested in $\cg$. Suppose now that there exist $i\in [t]$ so that $x_i\in B$. In this case we have that $x_j\leq \alpha$ for all $j\in [t]$, and $\alpha\geq \vee_{j=1}^t x_j>x_i$. Therefore $\vee_{j=1}^t x_j\in \cl\setminus \cg$. We conclude that  $B\cup  (A\setminus \{\alpha\})\in \cn(\cg)$.

$\bullet$ By definition, $f$ is injective.

$\bullet$ The map $f$ is surjective:

If $F\in \blb (\cf(\cn(\cg)))$, and $F$ is a face of $\cn(\cg)$, then $F\nsup B$. We want to prove that $F\in\cn(\cg')$. Let $\{y_1,\dots,y_s\}$, $s\geq 2$, be a set of incomparable elements in $F$. Since $\vee_{i=1}^s y_i\in \cl\setminus \cg$, and $\cg'=\cg\cup\{\alpha\}$, we are supposed to prove that $\vee_{i=1}^s y_i\neq \alpha$. Assume the contrary. Then, by \cite[Proposition 2.8]{FK}, $\{y_1,\cdots,y_s\}=\max \cg_{\leq y_1\vee\dots\vee y_s}=\max \cg_{\leq \alpha}=B$, which is a contradiction with $F\nsup B$. Hence $F\in \cf(\cn(\cg'))$, and $f(F)=F$.

If $[B,A]\in \blb(\cf(\cn(\cg)))$, then $A\nsup B$, and $A\cup B\in \cn(\cg)$. It is clear that $\alpha\notin A$, since $\alpha\notin \cg$. We ought to prove that $A\cup\{\alpha\}$ is nested in $\cg'$. Let us first prove that all elements in $A$ are comparable with $\alpha$. Suppose the contrary, i.e. that there exists\index{} $y\in A$, so that $y$ is not comparable with $\alpha$. Then $S:=B\setminus \cl_{\leq y}\neq \emptyset$, since otherwise $y\geq \bigvee B=\alpha$. Now, if $\alpha \vee y$ exists in $\cl$, we would have $\alpha<y\vee \alpha=y\vee \bigvee B=y\vee \bigvee S$, and hence $y\vee \bigvee S\in \cg$. This is a  contradiction with $A\cup B\in \cn(\cg)$, since $\{y\}\cup S\subseteq A\cup B$ is a set of incomparable elements of cardinality at least two.

Since any set of incomparable elements $\{a_1,\dots,a_s\}$ in $A\cup\{\alpha\}$, where $s\geq 2$, is actually a subset of $A$, by the same arguments as above we conclude that $\vee_{j=1}^s y_j\in \cl\setminus\cg'$. Therefore $A\cup\{\alpha\}$ is nested in $\cg'$ and $f(A\cup \{\alpha\})=[B,A]$.

$\bullet$ It is clear that $f$ is order-preserving.

Therefore, $f$ is an isomorphism of posets.
\end{proof}
\end{prop}
 From the previous two propositions, we can directly deduce the lattice case of \cite[Theorem 2.2]{BPSZ}:
 \begin{cor}
 Let $L$ be a lattice with finite length $l(L)<\infty$, and let $I\subset L$ be a proper ideal. Then $\Delta(\BBli)$ is obtained from $\Delta(\overline{L})$ by sequence of stellar subdivisions on all the edges from the set $\cs=\left\{\{x,y\}\mid x\in I\setminus\{\bm\},\ y\in \overline{L}\setminus I\right\}$, where these subdivisions are performed in an order of increasing length $l(x,y)$. 
 \begin{proof} If $\cg$ is a building set defined in Proposition \ref{prop_buildingBier}, then it is not hard to see that $\BBli=\cg\cup \cs'$, where $\cs'=\{[x,y]\mid \{x,y\}\in \cs,\textnormal{ with }x\in I\setminus \{\bm\}\textnormal{ and }y\in \overline{L}\setminus I\}$. Let $\cs=\{e_1,\dots,e_k\}$, where if
   $e_i=\{x_i,y_i\}$, $e_j=\{x_j,y_j\}$, and $i<j$, then
   $l(x_i,y_i)\leq l(x_j,y_j)$. Set $\cg_0=\cg$, and for $i\in [k]$,
   $\cg_i=\cg_{i-1}\cup \{[x_i,y_i]\}$, where $e_i=\{x_i,y_i\}$. It is
   clear that $[x_i,y_i]\in \max \left(\Blil\setminus
     \cg_{i-1}\right)$, and therefore, by Proposition
   \ref{prop_building_blowup},
   $\cf(\cn(\cg_i))=\bl_{\{[x_i,\tp],[\bm,y_i]\}}(\cf(\cn(\cg_{i-1})))$. In other words, since $\Delta(\overline{L})=\cn(\cg_0)$, and blowup in this case correspond to stellar subdivision of the edge $\{x_i,y_i\}$, we have that
 $$\cn(\cg_i)=\sd_{e_i}(\sd_{e_{i-1}}(\dots(\sd_{e_1}\Delta(\overline{L})))).$$
 We finish the proof remarking that $\cn(\Blil,\BBli)=\Delta(\BBli)$, since $\BBli$ is the maximal building set in $\Blil$.
 \end{proof}
\end{cor}
It is important to emphasize the following:
\begin{cor}\label{cor_PL}
For any lattice $L$, and a proper ideal $I\subset L$, $\Vert \Delta(\BBli)\Vert$ and $\Vert \Delta(\overline{L}) \Vert$ are PL homeomorphic. Furthermore, if $L$ is a face lattice of a strongly regular PL CW-sphere, then so is $\Bli$. 
\end{cor}

\section{Recursive coatom orderings and shellability of Bier lattices}\label{shell}

We now proceed to study the case of a shellable lattice or poset. After recalling the definition of shellability of a simplicial complex we will prove a proposition asserting that combinatorial blowups, and thus stellar subdivisions, preserve shellability.

\begin{DEF}\label{def_recursive}
Let $P$ be a bounded poset. We will say that $P$ admits a {\bf recursive coatom ordering} if $P=\{\bm,\tp\}$, or if there exists a coatom ordering $c_1,\dots,c_r$ so that the following two conditions are satisfied:
\begin{itemize}
\item[(R)] For all $j\in [r]$, the poset $[\bm,c_j]$ admits a recursive coatom ordering in which coatoms of  $[\bm,c_j]$ which are contained in $[\bm,c_i]$, for some $i<j$, come before all other coatoms.
\item[(S)] For all $1\leq i < k \leq r$ and all $x\in P$, if $x\leq c_i$ and $x\leq c_k$, then there exists some $j<k$ and some coatom $\omega$ of  $[\bm,c_k]$ so that $x\leq \omega \leq c_j$.
\end{itemize}
\end{DEF}
\begin{rem}
It was noticed in \cite[Proposition 2.13]{Sha} that, in the case when $L$ is a finite lattice and $c_1,\dots,c_r$
 is some coatom ordering of $L$, then this ordering satisfies condition (S) of Definition \ref{def_recursive} if and only if it satisfies the following condition:
 
\begin{itemize}
\item[(T)] {\em For all $1\leq i<j\leq r$, there exists some $k<j$ so that $$c_i\wedge c_j\leq c_k\wedge c_j\lessdot c_j.$$}
\end{itemize}
Since we will work with face lattices of simplicial complexes, we will verify conditions (R) and (T). The face lattice of a simplicial complex $K$ we will denote by $\hat{\cf}(K)$, that is $\hat{\cf}(K)=\cf(K)\cup\{\tp\}$.\end{rem}

We will now state two propositions which will be used further in the paper. For their proofs we refer the reader to the corresponding papers. 
 
\begin{prop}\tn{\cite[Theorem 5.13]{BWsn}}\label{shellRCO}
A simplicial complex $K$ is shellable \Iff  $\hat{\cf} (K)$ admits a recursive coatom ordering.
\end{prop}

\begin{prop}\tn{\cite[Theorem 5.1]{BWlsp}}\label{RCO}
A graded poset $P$ is totally semimodular if and only if for every interval $[x,y]$ of $P$, every atom ordering in $[x,y]$ is a recursive atom ordering.
\end{prop}

The next theorem is the main result of this section.

\begin{thm}\label{thm_rco}
Let $\hat{\cf}=\hat{\cf}(K)$ be a face lattice of some simplicial complex $K$, $\cf=\cfh\setminus\{\tp\}$, and assume that $\hat{\cf}$ admits a recursive coatom ordering. Then $\hat{\cf}_\alpha=\bl_\alpha ({\cf})\cup \{\tp\}$ also admits a recursive coatom ordering, where $\alpha$ is any element of $\hat{\cf}\setminus\{\bm,\tp\}$.
\begin{proof}

Since $\cfa$ is the face lattice of the simplicial complex $\sd_\alpha(K)$, intervals below maximal faces in $\cfa$ are Boolean. Having in mind that every boolean lattice is self-dual, graded, and totally semimodular (and that a recursive atom ordering of a poset is a recursive coatom ordering of its dual), by  Proposition \ref{RCO} every coatom ordering of these intervals is recursive. Therefore, to check that some ordering of coatoms of $\cfa$ is recursive, it suffices to prove that it satisfies condition  (T) from the remark above.

Suppose now that $\co=\{F_1,\dots,F_n\}$ is a \rco of $\cfh$.  Define $\cc \subset \cf$ in the following way: $$\mathcal{C}=\max\{G\in \cf\mid  G\nsupseteq \alpha\tn{ and }G\cup \alpha\tn{ is an element of }\cf\}.$$ Set $I=\{i\in [n]\mid F_i\supseteq \alpha\}$. Then it is not hard to see that all coatoms of $\cfa$ are $$\coat=\{F_i\mid i\in [n]\setminus I\}\cup\{[\alpha,G]\mid G\in \cc\}.$$ Denote the elements of $I$ with $i_1,\dots,i_t$, where $t=\vert I\vert$, and $i_1<i_2<\cdots<i_t$. Having in mind that the $F_i$'s are the coatoms of $\cfh$, and that $F_i\supseteq \alpha$, for all $i\in I$, it is not hard to see that $$\cc=\{G\in \cf\mid G\nsup \alpha\tn{ and }G\tn{ is a codimension }1\textnormal{ face of }F_i,\tn{ for some }i\in I\}.$$

Let now, for all $l\in [n]$, $\ca_l=\{G\in \cc\mid G\subseteq \overline{F_l}\setminus(\bigcup_{j=1}^{l-1}\overline{F_j})\}$. It is easy to see that $\cc=\bigsqcup_{1\leq l\leq i_t} \ca_l$. For any $G\in \cc$, let $A(G)=l$, where $G\in \ca_l$, and let $\mf(G)=\min\{s\in I\mid G\lessdot F_s\}$. It is clear that $A(G)\leq \mf(G)$.

Define a relation $\prec$ between coatoms of $\cfa$ in the following way:
\begin{itemize}
\item For any $i,j\in [n]\setminus I$, $F_i\prec F_j$ \Iff $i<j$.
\item For any $i\in [n]\setminus I$ and $G\in \cc$, $F_i\prec [\alpha,G]$ if and only if $i< \mf (G)$. 
\item For any $E,G\in \cc$, $[\alpha,E]\prec [\alpha,G]$ \Iff either $\mf(E)<\mf(G)$, or $\mf(E)=\mf(G)$ and $A(E)<A(G)$.
\end{itemize}
It is not hard to check that $\preceq$ is indeed a partial ordering. We will choose a linear extension thereof and denote it by $\coa$.  In order to prove that condition (T) holds for $\coa$, we need to prove the following four cases:

(1) Assume $F_i,F_j\in \coat$, and let $F_i$ come before $F_j$ in $\coa$. Then, by definition, $i<j$, and $\alpha\nsubseteq F_i,F_j$. Since $\co$ is a \rcok, there exists $k<j$ so that $F_i\cap F_j\subseteq F_k\cap F_j=F_j\setminus\{v\}$, where $v$ is some vertex of $F_j$.
\begin{enumerate}
\item[(1.1)] If $\alpha\nsub F_k$, then $F_k\in \coat$, clearly $F_k$ comes before $F_j$ in $\coa$, and $F_i\wedge F_j\leq F_k\wedge F_j\lessdot F_j$.
\item[(1.2)] If $\alpha\subseteq F_k$, pick any vertex $w\in \alpha\setminus F_j\subset F_k$. Then $F_k\setminus \{w\}$ is a codimension $1$ face of $F_k$ which does not contain $\alpha$, and therefore $[\alpha,F_k\setminus \{w\}]\in \coat$. Noticing that $\mf(F_k\setminus\{\alpha\})\leq k < j$, it is clear that $[\alpha,F_k\setminus \{w\}]$ comes before $F_j$ in $\coa$. Since $F_j\wedge [\alpha,F_k\setminus \{w\}]=F_j\cap (F_k\setminus \{w\})=F_j\setminus\{v\}$, we have that $F_i\wedge F_j\leq [\alpha,F_k\setminus \{w\}]\wedge F_j \lessdot F_j$.
\end{enumerate}

(2) Let $F_i,[\alpha,G]\in \coat$, and suppose that $F_i$ comes before $[\alpha,G]$ in $\coa$. Then $\alpha\nsub F_i$ and $i<\mf(G)$. Suppose first that $A(G)<\mf(G)$. Then $G\subseteq F_{A(G)}$, where $\alpha \nsub F_{A(G)}$, and $F_{A(G)}$ comes before $[\alpha,G]$ in $\coa$. Then we have that $F_i\wedge [\alpha,G]=F_i\cap G\leq  G=F_{A(G)}\cap G=F_{A(G)} \wedge [\alpha,G]\lessdot [\alpha,G]$, and condition (T) is satisfied in this case.

Let now $A(G)=\mf(G)=j$. By the assumption, there exists $k<j$ so that $F_i\cap F_j\subseteq F_k\cap F_j=F_j\setminus\{v\}$, for some $v\in F_j$. Note that $v\in G$, since otherwise $G\subseteq F_k$, and therefore $A(G)\leq k<j $, which would give a contradiction. 

Since $G\lessdot F_j$, denote with $w$ the vertex so that $F_j\setminus\{w\}=G$, where clearly $w\in \alpha$, and let 
\begin{equation}\label{eqn_H}
H=\left\{\begin{array}{ll} F_j\setminus\{v\}\lessdot F_j, & \tn{if }\alpha\nsub F_k; \\ F_k\setminus\{w\}\lessdot F_k, & \tn{if }\alpha\subseteq F_k. \end{array}\right.
\end{equation}

By simple checking, one can see that $\alpha\nsub H$, and $H\in \cc$, since $H$ is a coatom in some face that contains $\alpha$. Seeing that $H\subseteq F_k$, we conclude that either $\mf(H)<\mf(G)$, or $\mf(H)=\mf(G)$ and $A(H)<A(G)$, and therefore $[\alpha,H]$ comes before $[\alpha,G]$ in $\coa$. Now it is easy to see that $$F_i\wedge [\alpha,G]=F_i\cap G< G\setminus\{v\}< [\alpha,G\setminus\{v\}]=[\alpha,G]\wedge [\alpha, H]\lessdot [\alpha,G].$$

(3) If $[\alpha,G], F_j\in \coat$ and $[\alpha,G]$ comes before $F_j$ in $\coa$, then it must be $\mf(G)<j$. The proof for this case is the same as for case (1), having in mind that $ [\alpha,G]\wedge F_j=G\cap F_j\leq F_{\mf(G)}\cap F_j$, and since in case (1) we didn't use the fact that $\alpha \nsub F_i$.

(4) Finally, the case when $[\alpha,E],[\alpha,G]\in \coat$, and $[\alpha,E]$ comes before $[\alpha,G]$ in $\coa$. Then there are two possibilities, either $\mf(E)=\mf(G)$ or $\mf(E)<\mf(G)$.

Case when $j:=\mf(E)=\mf(G)$ is easy, namely both $E$ and $G$ are codimension $1$ faces in $F_j$, and hence $E\cap G$ is a codimension $1$ face of $G$. Therefore, we have that $[\alpha,E]\wedge [\alpha,G]=[\alpha,E\cap G]\lessdot [\alpha,G]$.

If $\mf(E)<\mf(G)$, let $i=\mf(E)$ and $j=\mf(G)$. Define $w$ to be the vertex so that $F_j\setminus \{w\}=G$, and $H$ as in equation (\ref{eqn_H}). Then, since $E\cap G\subseteq F_i\cap G\subseteq G\setminus\{v\}$, we have that 
\bdm [\alpha,E]\wedge [\alpha,G]=[\alpha,E\cap G]\leq [\alpha,G\setminus\{v\}]=[\alpha,G]\wedge [\alpha, H]\lessdot [\alpha,G].\edm

Therefore we have proved that $\coa$ satisfies condition (T) and hence is a recursive coatom ordering of $\cfa$. 
\end{proof}
\end{thm}

Since stellar subdivisions are described by combinatorial blowups in the face poset, we can formulate the following corollary.

\begin{cor}
If $K$ is a shellable simplicial complex, then so is $\sd_F(K)$, where $F$ is any nonempty face of $K$.
\end{cor}

Now we return to the general Bier poset construction, and conclude that:

\begin{cor}\label{cor_shellable}
Let $L$ be a lattice so that $\Delta(\overline{L})$ is shellable. If $I\subset L$ is a proper ideal in $L$, then $\Delta(\BBli)$ is also shellable. 
\end{cor}


\begin{Rem}{\bf(The general case)} {\em
If we consider any bounded poset $P$, the  first part of Corollary
\ref{cor_PL} remains true (see \cite[Corollaries 2.3 and 2.4]{BPSZ} for a
proof of this fact). Since stellar subdivisions are described by
combinatorial blowups in the face poset, Corollary \ref{cor_shellable} also remains true in the general case.

Let us also mention here that, by a remark at the end of Section 2 and  \cite[Theorem 2.2]{BPSZ}, for every bounded poset $P$ and for any proper ideal $I\subset P$, the simplicial complexes $\Delta(\overline{P})$ and $\Delta(\BBpi)$ have the same simple homotopy type.
}\end{Rem}

\begin{Rem}{\bf(Shellable spheres)} {\em It is clear that using the Bier poset construction together with
Corollaries \ref{cor_PL} and \ref{cor_shellable}, one can obtain
numerous simplicial shellable $n$-spheres with more than $2(n+2)$
vertices, therefore answering  one of the questions asked in
\cite{BPSZ}. 

Choosing an appropriate poset $L$ we also obtain numerous shellable
simplicial  complexes with the homotopy type of wedges of spheres, with any number of vertices.  
}\end{Rem}


\end{document}